\documentclass[10pt]{amsart}
\usepackage{epsf,amssymb}

\newtheorem{theorem}{Theorem}

\theoremstyle{remark}

\renewcommand{\leadsto}{\rightarrow}
\newcommand{\vpic}[1]{{\rule[-.7em]{0em}{2em}\:\vcenter{\epsffile{#1.eps}}\:}} 
\newcommand{\thetadiag}{\vcenter{\epsffile{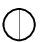}}}
\newcommand{\str}{\rule[-.7em]{0em}{2em}}

\newcommand{\Q}{\mathbb{Q}\,}
\newcommand{\R}{\mathbb{R}\,}
\newcommand{\C}{\mathbb{C}\,}
\newcommand{\A}{\mathcal{A}}
\newcommand{\Z}{\mathbb{Z}}
\newcommand{\f}{{{}^f\!}}

\newcommand{\fA}{{}^f\!\!\mathcal{A}}

\begin{document}
\thispagestyle{empty}
\vbox{\vskip 20pt}
{\center {{\bf THE KONTSEVICH INTEGRAL AND ALGEBRAIC\\ STRUCTURES ON THE
SPACE OF DIAGRAMS}
\\
\bigskip
{\small {SIMON WILLERTON}%
\\
\smallskip 
{\em  Institut de Recherche Mathématique Avancée, Université
Louis Pasteur et CNRS,\\ 
7 rue René Descartes, 67084 Strasbourg, France}
\vskip 1in
\setlength{\baselineskip}{10pt}
\footnotesize\begin{minipage}{4.5in}\section*{abstract}
This paper is part expository and part presentation of calculational
results.  The target space of the Kontsevich integral for knots is a
space of diagrams; this space has various algebraic structures which
are described here.  These are utilized with Le's theorem on the
behaviour of the Kontsevich integral under cabling and with the
Melvin-Morton Theorem, to obtain, in the Kontsevich integral for torus
knots, both an explicit expression up to degree five and the general
coefficients of the wheel diagrams.
\\{\em Keywords:\/} Vassiliev invariants, torus knots, Kontsevich integral.
\end{minipage}}}\\}
\medskip
\section*{Introduction}
The Kontsevich integral, $Z$, is a map from the set of
oriented knots to a space, $\A$, of ``diagrams''.  For example it will
be seen in Section~8 that the value taken on a negative trefoil ---
the $(2,-3)$-torus knot --- begins
$$Z\left(\vpic{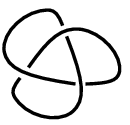}\right)=
\vpic{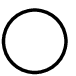} - \tfrac{1}{2}\vpic{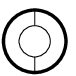} -\tfrac{1}{2}\vpic{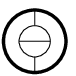}
-\tfrac{31}{48}\vpic{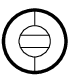} +\tfrac{5}{24}\vpic{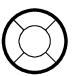}+
\tfrac{1}{8}\vpic{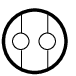}+\dots.$$ 
The Kontsevich integral was introduced in
\cite{Bar:OnVIs,Kontsevich:Vassiliev}; it is universal for
Vassiliev knot invariants in the sense that every finite-type
invariant factors through the Kontsevich integral, and it is universal
in a certain sense
for knot invariants coming from quantum groups.  Other approaches to
defining this invariant have been introduced 
but none give in easily to attempts at actual calculation.

The space of
diagrams has a rich algebraic structure, the following will be
described in this paper:  it has two
products making it a graded Hopf algebra in two ways, it has a set of
simultaneously diagonalizable ``Adams operations'' and the space of
primitive elements in $\A$ has the structure of a module over Vogel's
algebra $\Lambda$.  How all of these structures interact and how they relate to
topological properties of knots is not yet properly understood.
However, they can be utilized to describe the Kontsevich integral for
torus knots and to identify the Alexander polynomial in terms of wheel
diagrams.  These last two things are taken in this paper as
motivations for the algebraic structures described herein. 

The paper is organized as follows.  
The first section contains a brief sketch of the original definition
of the Kontsevich integral.  The second section introduces the target
space $\A$ and gives an explicit basis up to degree five --- the
specific choice of basis is made to interact nicely with the algebraic
structure described later.  The Hopf algebra structure of $\A$ coming
from the connect sum product is presented in Section~3, together with
some consequences of the Hopf algebra structure theorems, such as
identifying $\A$ as being freely generated by the space spanned by
connected diagrams.  Section~4 shows how the theory naturally extends
to framed knots and how the target space is extended to
$\fA$ by the addition of an extra generator to encode the framing.  In
Section~5 the Adams operations are introduced, as are symmetrized
diagrams which are eigenvectors for these operations; also wheel
diagrams are defined and used with the disjoint union product to give
a description of the element $\Omega$.  The module structure of the
primitive elements of $\A$ over
Vogel's algebra $\Lambda$ is the subject of the sixth section.
Cabling knots and Le's theorem on the behaviour of the Kontesevich
integral are contained in the seventh section.  Section~8 presents the
result of calculating the Kontsevich integral of torus knots by using
the machinery so far introduced, and also gives a criterion for
Vassiliev invariants of torus knots.  The final section shows how the
Melvin-Morton Theorem allows the identification of the wheels
coefficients from the Alexander polynomial and, in particular, does
this explicitly for torus knots.  Appendix~A lists the change of basis
data from round diagrams to symmetrized diagrams used in Section~8 and
Appendix~B compares the results here with those of Alvarez and Labastida.

This paper is an extended version of the talk that I gave at Knots in
Hellas '98, and this supersedes the earlier preprint
\cite{Willerton:CablingTorus}. 

\section{The Kontsevich integral.}
Here I will try to give a flavour of the original definition of the
Kontsevich integral --- it still remains somewhat mysterious --- one of the
goals of this paper is to show that it can be calculated for certain
knots without performing any sorts of integral.  The target space
$\A$ will not be defined until the next section, so here it should be
thought of as ``linear combinations of diagrams''.

Suppose that $K$ is a Morse knot, that is a 
knot embedded in $\R^3=\R_t\times\C_z$ such that the projection to the $\R$
co-ordinate is a Morse function. 
\begin{figure}
$$\vpic{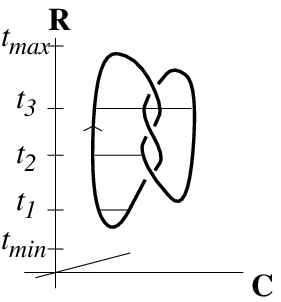}\quad\leadsto\quad D_P=\vpic{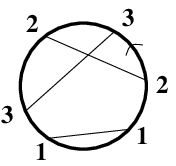}$$
\caption{A Morse knot together with a set of pairings $P$ for $m=3$,
and the corresponding diagram $D_P$.}
\label{KontIntFig}
\end{figure}
Then define
$$\overline z (K):=\sum_{m=0}^\infty \frac{1}{(2\pi i)^m}
\int\limits_{\genfrac{}{}{0pt}{1}{t_{\text{min}}<  t_1 < \dots <
t_m< t_{\text{max}}}{t_j \text{ non-critical}}}
\sum_{\genfrac{}{}{0pt}{1}{\text{pairings }}{P=(\{z_j,z_j'\})}}(-1)^{\#\downarrow P}
[D_P]\bigwedge_{j=1}^m \frac{dz_j-dz'_j}{z_j-z'_j}.$$
where a set of pairings $P$ is a set of unordered pairs
$\{z_j,z_j'\}$ so that 
$(t_j,z_j),(t_j,z_j')\in \R\times \C$ are distinct points on the knot
on the same level,  where $\#\!\downarrow\! P$ is the number of points
in the pairing such that the orientation is pointing downwards and
where $D_P$ is the diagram obtained by drawing the preimage of the
knot --- a circle --- and drawing chords connecting the points
corresponding to the pairs.  See Figure~\ref{KontIntFig}.  

Considered as an element of $\A$, this gives an invariant of Morse
knots, but is not invariant if one 
introduces a ``hump'', i.e.\ a minima-maxima pair with respect to the
$\R$-axis, which increases the number of critical points.  
It is necessary to introduce the following correction
term:  define $\Omega:=\overline
z\left({\vcenter{\epsffile{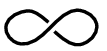}}} \right)^{-1}\in \A$,
where this is a Morse embedding of the unknot with four critical
points.  $\Omega$
is an important element which will also 
be returned to later in the paper.  Now the
Kontsevich integral for a Morse knot $K$ with $c(K)$ critical points
is defined to be
$$Z(K):=\overline z(K).\Omega^{c(K)/2-1}\in\A.$$
For a pleasant survey going in to some of
the finer points of this definition, 
the reader is directed \cite{ChmutovDuzhin:KontsevichIntegral}. 

It is important to note here that there are two normalizations
of the Kontsevich integral in common use: the other can be defined by 
$Z_q(K):=\overline z(K).\Omega^{c(K)/2}$.  It is clear from this that
the two normalizations differ by a factor of $\Omega$ and that
$Z(\text{unknot})=\vcenter{\epsffile{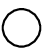}}$, while $Z_q(\text{unknot})=\Omega$, and so
$\Omega$ is often referred to as ``the Kontsevich integral of the
unknot''.  Both normalizations have their advantages: $Z$ is
multiplicative under connect sum, while $Z_q$ is better behaved under
cabling operations.  In this work I will use exclusively the
multiplicative normalization.

There are a couple of other ways of defining such a knot invariant
(see \cite{BarStoi:Odense}):
one is via quasi-Hopf algebras and Drinfeld associators, this is
essentially a discretized version of the above; the other way is the
so-called Bott-Taubes approach using configuration space integrals,
this is believed to give the same answer as the Kontsevich
integral above.

\section{The space of diagrams $\A$.}
The target space $\A$ of the Kontsevich integral will be defined here,
and a basis given up to degree five.

\begin{figure}[t]
$$\vpic{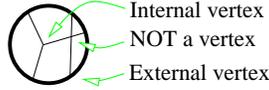}$$
\caption{A degree three diagram.  Note that the four-valent vertex is
not a vertex of the diagram, but just a result of the planar drawing.}
\label{diagramfig}
\end{figure}

First it is
necessary to define the notion of a {\em diagram\/}  (which I will also call a
round diagram) --- this was called a Chinese character diagram in
\cite{Bar:OnVIs}.  A diagram is a trivalent 
connected graph such that each vertex is cyclically oriented and such
that there is a distinguished oriented cycle called the {\em Wilson loop\/}.
A diagram is graded by half the number of vertices that it has.  By
convention a diagram is drawn such that the Wilson loop is an
anticlockwise oriented circle, and such that each vertex is also
oriented anticlockwise.  See Figure~\ref{diagramfig}.  The complement
of Wilson loop is termed the {\em internal graph}.

Define the rational vector space $\A_n$ to be the vector space
spanned by the diagrams of degree $n$ and subject to the so-called
STU and 1T relations:
$$\A_n:=\langle \text{diagrams of degree\ }n\rangle_\Q /\left\{
\begin{array}{rcll}\vpic{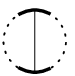}&=&0:&\text{1T}\\
\vpic{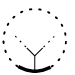}&=&\vpic{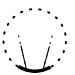}-\vpic{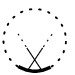}:&\text{STU}.\end{array}\right.$$
Note that the 1T relation means that a diagram vanishes if it has a
chord which is not intersected by any other part of the internal
graph.  It is also worth noting that the antisymmetry (AS) relation is a direct
consequence of the STU relation: antisymmetry says that reversing the
orientation of one vertex in a diagram is the same as multiplying the
diagram by minus one; pictorially,
$$\vpic{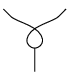}=-\vpic{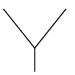}:\quad \text{AS.}$$

Look now at the spaces of low order, and find explicit bases.
Degree zero is very easy:
$$\A_0:=\left\langle \vpic{1A}\right\rangle_\Q.$$
Degree one is almost as easy:
$$\A_1:=\left\langle \vpic{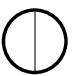}\right\rangle_\Q/\left 
                      (\vpic{IA}=0\right)\cong 0.$$
It starts getting interesting at degree two:
$$\A_2:=\left\langle \vpic{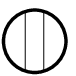},\vpic{w2A},\vpic{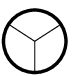},\vpic{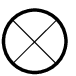}
\right\rangle_\Q/ \left\{
\begin{array}{l} \vpic{I2A}=0\\\vpic{tree3A}= \vpic{I2A}-\vpic{XA}\\
\vpic{w2A}=2\vpic{tree3A}.\end{array}\right.$$
That has four generators and three relations, so pick the following
generator:
$$\A_2\cong \left\langle\vpic{w2A}\right\rangle_\Q.$$
Degree three is more awkward:
$$\A_3:=\left\langle
\vpic{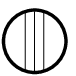},\vpic{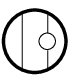}, \vpic{tw2A}, \vpic{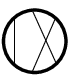}, \vpic{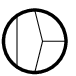},
\vpic{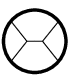}, \vpic{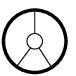}, \vpic{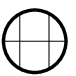}, \vpic{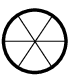}
 \right\rangle_\Q/\text{relations}.$$
It transpires that $\A_3$ is also one dimensional and the generator
can be chosen to be as follows:
$$\A_3\cong \left\langle\vpic{tw2A}\right\rangle_\Q.$$
For later use, we will just write down bases for $\A_4$ and $\A_5$:
$$\A_4\cong \left\langle\vpic{t2w2A},\vpic{w4A},\vpic{w2w2A}
     \right\rangle_\Q;\quad
\A_5\cong
\left\langle\vpic{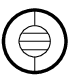},\vpic{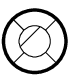},\vpic{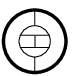},\vpic{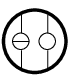}
      \right\rangle_\Q.$$ 
A rationale for this choice of basis should become apparent during the
of course of the paper; it
has primarily been chosen to simplify the calculation of the results
in Appendix~A.  
The dimension of $\A_n$ is only known for $n\le 12$ and these
dimensions are tabulated in Table~\ref{dimtable}.  See
\cite{Dasbach:CombinatorialStructureIII} for some results on asymptotic lower
bounds on the dimensions of these spaces.
\begin{table}
\begin{tabular}{|r|c|c|c|c|c|c|c|c|c|c|c|c|c|}
\hline $n$&0&1&2&3&4&5&6&7&8&9&10&11&12\\
\hline
$\dim \A_n$ &1&0&1&1&3&4&9&14&27&44&80&132&232\\
\hline
\end{tabular}
\smallskip \caption{The dimensions of the spaces $\A_n$ up to $n=12$
as given in \cite{Bar:OnVIs, Kneissler:VasDegreeTwelve}.} 
\label{dimtable}
\end{table}

Define $\A$ to be the projective limit of the direct sum 
$\bigoplus_{n=0}^\infty \A_n$, 
so an element of $\A$ is a possibly infinite linear combination
of diagrams which is finite in each degree.  This space has lots of
algebraic structure which will be described in the sections below.

\section{Hopf algebraic structure of $\A$.}
In this section I will introduce a Hopf algebra
structure for $\A$ coming from the connect sum product,
and show consequences of the Hopf algebra structure
theorems, notably that the Kontsevich integral of a knot can be
expressed as an exponential of a sum of connected diagrams. 
\subsection{Product and coproduct on $\A$.} 
Firstly, $\A$ has a natural product $\A\otimes\A\to\A$.  
This is defined on diagrams in
terms of a connect sum operation.\\ E.g.
$$\vpic{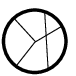}.\vpic{XA}=\vpic{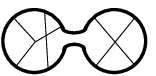}=\vpic{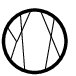}
$$
The resulting diagram is, of course, dependent on where the cuts are
made in the two diagrams; however, the result is well-defined modulo
the STU relation (see \cite{Bar:OnVIs}).  This makes $\A$ into a
graded, commutative algebra.%
\footnote{This means commutative in the genuine sense and not the
commutative-graded sense.  Philosophically this is because $\A$ all
lives in homological degree zero.  The structure theorems of
\cite{MilnorMoore} deal with commutative-graded algebras, but this is
not a problem as, if necessary, the grading of $\A$ can be doubled so
that everything lives in even degree.}
The unit is the diagram which has no internal graph.

This product is related to the connect sum operation, $\#$, on knots via the
Kontsevich integral, as the following theorem demonstrates.
\begin{theorem}[Kontsevich]
For $K$ and $L$ knots:
$Z(K\# L)=Z(K).Z(L)$.
\label{MultiplicativityThm}
\end{theorem}

Secondly, $\A$ has a coproduct $\Delta\colon\A\to\A\otimes\A$. This is
defined on a diagram by taking all of the ways that the internal graph
can be decomposed into an ordered pair of complementary sets of
components.\\ E.g.
$$\Delta\left(\vpic{KA}\right)= \vpic{KA}\!\otimes\!\vpic{1A}+
   \vpic{tree3A}\!\otimes\!\vpic{IA}+\vpic{IA}\!\otimes\!\vpic{tree3A} 
  +\vpic{1A}\!\otimes\!\vpic{KA}.
$$
The counit is the function $\A\to\Q$ which maps the diagram with no
internal graph to $1$ and which vanishes on all other diagrams.  This
coproduct is graded and cocommutative, and is also compatible with the
product in the sense that it is a map of algebras:
$\Delta(a.b)=\Delta(a).\Delta(b)$.  This makes $\A$ into a connected,
graded Hopf
algebra {\it à la\/} Milnor and Moore%
\footnote{
This implies that it is uniquely equipped with an antipode, but this
is not really of concern here.}
\cite{MilnorMoore}.

\subsection{Consequences of Hopf algebra structure theorems.}
Now we can utilize the nice structure theorems of commutative,
cocommutative Hopf algebras.  

There are two special subsets of a Hopf
algebra.  The first is the vector space
of {\em primitive elements\/}; this is defined to be the set of
elements which satisfy $\Delta(x)=x\otimes 1 + 1\otimes x$.  This
space naturally forms a Lie algebra under the commutator bracket, but
in the commutative case of $\A$ the bracket is trivial.  The space of
primitive elements of $\A$ certainly contains those
 diagrams whose internal graph is connected --- the
so-called {\em connected\/} diagrams.
$$\text{Connected diagrams:\ }\vpic{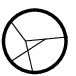},\vpic{tree4A}.\qquad
\text{Non-connected diagrams:\ }\vpic{tw2w2A},\vpic{KA}.$$
In fact by the following theorem it suffices to know these diagrams.
\begin{theorem}[Bar-Natan]
The vector space of primitive elements of $\A$ is (projectively) spanned
by connected diagrams.
\end{theorem}
One of the key structure theorems of Milnor and Moore is that a
commutative and cocommutative Hopf algebra is (naturally isomorphic
to) the polynomial algebra over its primitive elements.  Translating
this into the current context we obtain:
\begin{theorem}
The algebra $\A$ is the completed polynomial algebra on the
connected diagrams.
\end{theorem}
In particular this means that if one has a basis for the space of
connected diagrams then this can be extended to a basis of $\A$ by
taking the monomials in the basis of connected diagrams.  Now look
back to the bases given in Section~2 and see that this is precisely
what I have done.  The basis elements are all connected except one in
degree four and one in degree five which are products of lower degree
connected basis elements.

The second important subset of a Hopf algebra is the set of
{\em group-like elements\/}; these are the elements which satisfy
$\Delta(g)=g\otimes g$.  This set forms a group under the product; the
inverse operation is just the antipode map.  This set is important in
the current context because of the following theorem.
\begin{theorem}[Kontsevich]
The Kontsevich integral maps knots to group-like elements.
\label{ThmGrpLike}\end{theorem}

A theorem of Quillen \cite{Quillen:RationalHomotopy}
asserts that, under reasonable hypotheses, the usual exponential map,
$\exp(x):=\sum_{i=0}^\infty x^i/i!$, defines a bijection of sets
from the primitive elements to the group-like elements, with the
inverse being given by logarithm.
Combining with Theorem~\ref{ThmGrpLike} this gives
\begin{theorem}
The Kontsevich integral of a knot is expressible as the exponential of
a sum of connected diagrams.
\end{theorem}
\noindent For instance the example from the introduction can be rewritten as 
$$Z\left(\vpic{smalltrefoil}\right)
=\exp_{\mbox{.}}\!\left(-\tfrac{1}{2}\vpic{w2A}-\tfrac{1}{2}\vpic{tw2A}
-\tfrac{31}{48}\vpic{t2w2A}+\tfrac{5}{24}\vpic{w4A}+ \dots\right).$$
Here $\exp_{\mbox{.}}$ means that the connect sum product is used in the
definition of the exponential map; this is to prevent confusion when
the disjoint union product, $\amalg$, is introduced later.

The Baker-Campbell-Hausdorff formula describes how the product on the
group-like elements is pulled back to the primitive elements; in the
commutative case this reduces to just addition on the space of
primitive elements.  Together with Theorem~\ref{MultiplicativityThm},
this means that not only is 
$\ln_{\mbox{.}}\!\circ Z$ well defined on knots, but that it is additive under
connect sum of knots:
$\ln_{\mbox{.}}\!\circ Z(K\# L)=
   \ln_{\mbox{.}}\!\circ Z(K)+\ln_{\mbox{.}}\!\circ Z(L)$.

\section{Framed version of the theory.}
The Kontsevich integral can be naturally defined for framed knots.  The
target space of this has one extra generator which lives in degree
one --- this keeps count of the framing.
\subsection{Framed knots.}
\begin{figure}
$$\vpic{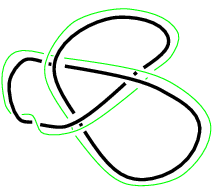}\leadsto\vpic{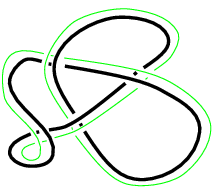}\leadsto\vpic{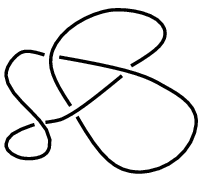}$$
\caption{This framed trefoil can be drawn in the blackboard framing by
adding a kink.}
\label{blackboardfig}
\end{figure}

A framed knot is a knot equipped with a homotopy class of sections of
the normal bundle.  This can be thought of as a parallel strand
running very close to the knot.  A knot diagram is said to be
blackboard framed if the framing curve lies parallel to the knot in
the plane of the paper.  Any framed knot can be drawn in the
blackboard framing, possibly by adding kinks, see
Figure~\ref{blackboardfig}.

For a given knot, the framings are classified by the framing number
(also called the writhe) which is the linking number of the knot and
its framing curve: for a framed knot $K$, this will be denoted by
$F(K)$.  

\subsection{The Hopf algebra $\fA$.}
The Kontsevich integral for framed knots will be defined below to take
values in a Hopf algebra of diagrams defined similarly to $\A$, but
without the 1T relation, namely:
$$\fA:=\langle\text{diagrams}\rangle_{\Q}/\text{STU}.$$
Thus $\A\cong\A/\text{1T}$.  In fact, $\A$ can be identified as a
sub-Hopf algebra of $\fA$ by an inclusion
$\varphi\colon\A\hookrightarrow\fA$ \cite{Willerton:Hopf}.  
Furthermore $\fA$ has just one
extra primitive generator --- this is the diagram in degree one:
$\thetadiag$.  So there is the natural isomorphism
$\fA\cong\A\oplus\thetadiag.\fA$.

\subsection{The framed Kontsevich integral.}   
The Kontsevich integral for oriented framed knots, taking values in $\fA$,
can be defined%
\footnote{Note that this will not quite give a universal Vassiliev invariant
of framed knots, a term involving the parity of the framing is
required --- see \cite{KasselTuraev}.}
by
$$\f Z(K):= Z(K).\exp_{\mbox{.}}\!\left(
           F(K)\thetadiag/2\right).$$ 
This looks a little {\it ad hoc}, but it is equivalent to other
definitions, see \cite{LeMurakami:Framed}.  Note that this extra
generator of $\fA$ is encodes precisely 
the framing data and on  the level of logarithms
$$\ln_{\mbox{.}}\!\circ\f Z(K)=\ln_{\mbox{.}}\!\circ Z(K)+F(K)\thetadiag/2.$$

\section{Adams operations and symmetrized diagrams.}
\subsection{Adams operations.}
On the algebra $\fA$ there is a set of operations, $\{\psi^m\}_{m\in
\Z}$, called the Adams operations.  For $m$ a
positive 
integer, $\psi^m$ is defined on a diagram by lifting the vertices on
the Wilson loop to the $p$th connected cover of the Wilson loop.\\E.g.
$$\psi^3\left(\vpic{XA}\right)=\text{``}\vpic{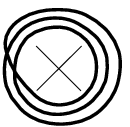}\text{''}
           =33\vpic{XA}+48\vpic{I2A}.$$ 
For $m$ a negative integer, one takes $\psi^{-m}$ of the diagram with
the orientation of the Wilson loop reversed --- note it is unknown
whether reversing the orientation of the Wilson loop induces the identity
map on $\fA$ or not.

The Adams operations satisfy the relation
$\psi^m\circ\psi^p=\psi^{mp}$, but they are not algebra maps, so in
general $\psi^m(a.b)\ne \psi^m(a)\psi^m(b)$.  However the  $\{\psi^m\}$
are co-algebra maps and in fact $\fA$ can be equipped with a different
multiplication for which the $\{\psi^m\}$ are Hopf algebra maps  ---
see below.

\subsection{Symmetrized diagrams and the disjoint union product.}
In certain situations, such as doing calculations involving the Adams
operations or for expressing the element $\Omega$, it is convenient to use a basis of $\A$ which is an
eigenbasis 
for the Adams operations; such a basis can be given in terms of
{\em symmetrized diagrams\/} (in \cite{Bar:OnVIs} they
are called Chinese characters). 

 A symmetrized diagram  is a graph with univalent
vertices, and cyclicly
ordered, trivalent vertices, such that each
connected component of the graph has at least one univalent vertex.
The univalent vertices are also called {legs\/}.  The element of 
$\fA$ that such a graph represents is obtained by averaging all
possible ways of attaching the univalent vertices to a Wilson loop.\\
E.g. 
$$\vpic{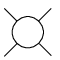}=\tfrac{1}{24}\left( 8 \vpic{w4A}+16\vpic{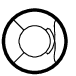}\right).$$

It is not difficult to show that these are eigenvectors for the Adams
operations: a symmetrized diagram with $u$ univalent vertices, is an
eigenvector for the Adams operation $\psi^p$ with eigenvalue $p^u$.
Also it is not difficult to show that symmetrized diagrams span $\fA$:
one can use the easy fact that a diagram can be written as the sum of
the symmetrized diagram obtained by removing the Wilson loop and
symmetrized diagrams with fewer legs.  If one defines the following IHX
relation:
$$\vpic{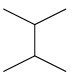}=\vpic{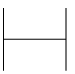}-\vpic{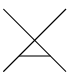}: \text{ IHX},$$
then one can in fact take the 
vector space generated by
symmetrized diagrams modulo the IHX and AS relations, and the above map
to $\fA$ induces an isomorphism.
Note that the connected symmetrized diagrams span the
primitive subspace of $\fA$.

The disjoint union of symmetrized diagrams, $\amalg$, 
gives $\fA$ a {\em second product\/} which is also
compatible with the co-product; thus it gives $\fA$ a second Hopf
algebra structure.   With respect to this second structure the
Adams operations are Hopf-algebra maps.  There seems not to be any
topological interpretation of this disjoint union product; however, it
can be related to the connect sum product by the (complicated) {\em
wheeling\/} map \cite{BarLeThurston:BLT}.

\subsection{Wheels.}
Certain symmetrized diagrams --- the {\em wheels\/} --- play a
particularly important r\^ole in the theory, both in the Alexander
polynomial (see Section~9) and the explicit formula for the element
$\Omega$ described below. 

If one considers the connected symmetrized
diagrams of degree $n\ge 2$ (i.e.\ with a total of $2n$ vertices) then
those with more than $n$ univalent vertices are trivial in $\fA$.  If
$n$ is odd then those with $n$ univalent vertices also vanish; but if
$n$ is even then define the non-trivial diagram  $w_n$ to be
the {\em wheel with $n$ legs}.
\\E.g.
$$w_2:=\vpic{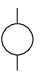};\qquad w_4:=\vpic{w4B}.$$
Each $w_n$ spans the ``top'' eigenspace of the degree $n$ primitive
diagrams for the Adams operations.

These wheels can be used together with the disjoint union product to
give a closed formula for the element $\Omega\in \A$ defined in
Section~1.  Firstly define the modified Bernoulli numbers,
$\{b_{2n}\}_{n=1}^\infty$  via
$$\sum_{n=1}^\infty b_{2n}x^{2n} 
   = \tfrac{1}{2} \ln \frac{\sinh (x/2)}{x/2}.$$
The following was conjectured in \cite{BarGar:Wheels} and proved in
\cite{BarLeThurston:BLT}:
\begin{theorem}[Bar-Natan, Le, Thurston]
$\Omega=\exp_\amalg \left(\sum\nolimits_{n=1}^\infty b_{2n} w_{2n}\right)$.
\end{theorem}
\noindent For instance, writing $1$ for the empty symmetrized
diagram of degree zero, 
 $$\Omega=1+\tfrac{1}{48}\vpic{w2B}-\tfrac{1}{5760}\vpic{w4B}
 +\tfrac{1}{2.48^2}\vpic{w2B}\vpic{w2B}+\{\text{terms of degree }\ge 6\}.$$
  This will be needed in the calculation in Section~8.

\section{Primitive elements as a module over Vogel's algebra $\Lambda$.}
The space of primitive elements of $\A$ admits a module structure,
introduced in \cite{Vogel:AlgebraicStructures}, which
turns out to  be very useful for performing low order calculations by
hand; it is also very useful in the context of weight systems
coming from Lie algebras, and it was used by Vogel to prove that not
all weight systems come from simple Lie super-algebras ---  this
will not be addressed further here, but the reader is directed to
\cite{Lieberum:VIsNotFromLAs, Kneissler:VasDegreeTwelve,
Kneissler:Relations, Vogel:AlgebraicStructures}.

It is important to note that I am working over the rationals, so I
denote by $\Lambda$ what is denoted $\Lambda\otimes\Q$ in
\cite{Vogel:AlgebraicStructures}.

Define a Vogel diagram to be a graph with cyclicly oriented trivalent
vertices  and three univalent vertices labelled bijectively with the
set $\{1,2,3\}$.  A Vogel diagram has degree $(\#\text{trivalent
vertices}-1)/2$.   Consider the space of linear combinations of
Vogel diagrams modulo the
IHX and AS relations, then $\Lambda$ is defined to be the subspace
consisting of all elements $u$ 
such that%
\footnote{If one is working over a ring in which $2$ is not invertible
then an extra condition is required on the elements of $\Lambda$, see 
\cite{Vogel:AlgebraicStructures}.}
if $\sigma$ is a permutation of the label set then
$\sigma(u)=\text{signature}(\sigma)u$.

Given a Vogel diagram and a connected diagram in $\A$, one can form a
new connected diagram by putting the Vogel diagram in the place of a
trivalent vertex, taking note of orientations.  For example
$$\vpic{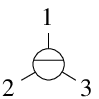},\vpic{tree4A}\qquad\leadsto\qquad \vpic{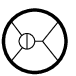}.$$
The constraint on the elements of $\Lambda$ ensures precisely that
this descends to a well-defined action of $\Lambda$ on the primitive
elements of $\A$, for instance it does not depend on which trivalent
vertex is replaced.  Similarly $\Lambda$ acts on itself, giving it a graded
algebra structure.

Up to degree six, the algebra $\Lambda$ is freely generated by the following
elements:
$$\tfrac{1}{2} x_1:=t:=\vpic{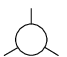};\quad x_3:=\vpic{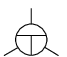};\quad x_5:=\vpic{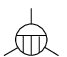}.$$
It was conjectured by Vogel that $\Lambda$ is freely generated by
$\{x_{2i-1}\}_{i=1}^\infty$ but Kneissler
\cite{Kneissler:VasDegreeTwelve} proved a relation between them
in degree ten: it is unknown whether these elements generate $\Lambda$.

Kneissler also showed that in degrees less than six, the primitive space of
$\A$ is freely generated over $\Lambda$ by the diagrams
$\vcenter{\epsffile{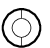}}$ 
and $\vcenter{\epsffile{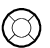}}$;  although more than just the
round wheel diagrams are 
needed in degree six.  The reader can check that this is precisely the
form that my basis of primitive elements takes.

I do not know how to sensibly extend this to a $\Lambda$ action on the
whole of $\A$, nor how to interpret the action topologically.

\section{Cabling knots and the theorem of Le}
Cabling is a natural operation on framed knots, with respect to which
the Kontsevich
integral is particularly well behaved, as is demonstrated
by the theorem of Le described below.

\begin{figure}[t]
$$\text{(i)}\ T_m=\vpic{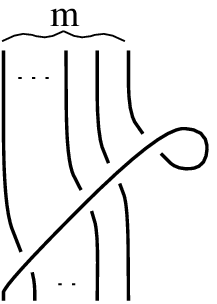};
\qquad\text{(ii)}\vpic{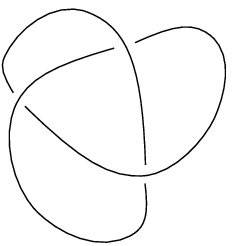}\stackrel{(3,2)}{\mapsto} \vpic{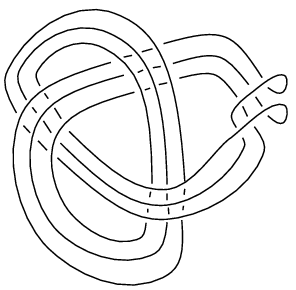}
$$
\caption{(i) The tangle used in the definition of the $(m,p)$-cable.
(ii) The $(3,2)$-cable of a blackboard framed trefoil.}
\label{CableFigure}
\end{figure}

For  $m$ and $p$ coprime integers and $K$ a framed knot, the
$(m,p)$-cable, $K^{(m,p)}$, of $K$ can be defined by taking a
blackboard framed knot diagram of $K$, drawing $m$ parallels of the
knot in the plane of the paper and inserting $p$ copies of the tangle
$T_m$ shown in Figure~\ref{CableFigure}(i), to obtain a blackboard framed
knot diagram representing $K^{(m,p)}$ --- see Figure~\ref{CableFigure}(ii).

The algebraic structure of $\fA$ described in the previous sections,
together with the element $\Omega$, permits
the following description of how the framed Kontsevich integral
behaves under cabling: 
\begin{theorem}[Le \cite{Le:Unpublished}]
For $m$ and $p$ coprime integers, and $K$ a framed knot,
$$\f Z(K^{(m,p)}).\Omega = 
\psi^m\left(\f Z(K).\Omega.
\exp_{\mbox{.}}\!\left(\tfrac{p}{m}\tfrac{1}{2}\thetadiag\right)\right).$$ 
\end{theorem}

Note that if  any of the weight system coming from 
a semi-simple Lie algebras is applied to both sides, then the
consequent equality follows from an analysis of the theorem in
\cite{Morton:ColouredJonesTorus} which Morton attributes to Strickland
\cite{Strickland:Unpublished} and Rosso and Jones
\cite{RossoJones:Torus}, and it was because of this that I conjectured
this theorem in an earlier version of this note before discovering
that Le had already proved it.

This theorem will be used in the next section to calculate the 
Kontsevich integral up to degree five for torus knots.
 
\section{The Kontsevich integral and Vassiliev invariants of torus knots.}
If $m$ and $p$ are coprime integers then the $(m,p)$-torus knot,
denoted $T(m,p)$, is the knot type that has a representative which is
embedded on the standard torus, wrapping $m$ times around
longitudinally and $p$ times around meridinally.  For instance, the
positive trefoil is the $(2,3)$-trefoil.  The $(m,p)$-torus knot
can equivalently
be described as the $(m,p)$-cable of the zero-framed unknot.  This
latter description allows the calculation of the Kontsevich integral
of the $(m,p)$-torus knot via the previous theorem, as the framed
Kontsevich integral of the zero-framed unknot is just the diagram in
degree zero.

To calculate $Z(T(m,p))$ proceed as follows: convert the expression
for $\Omega$ given in Section~5 to round diagrams using the change of
basis data in Appendix~A; form the connect sum product with
$\exp_{\mbox{.}}\!\left(\tfrac{p}{m}\tfrac{1}{2}\thetadiag\right)$; convert back to
symmetrized diagrams; apply $\psi^m$, recalling that symmetrized
diagrams are eigenvectors; convert back to round diagrams and form the
connect sum product with $\Omega^{-1}$; finally change to the zero
framing --- this will require multiplying by
$\exp_{\mbox{.}}\!\left(-mp\:\thetadiag/2\right)$. 

This calculation is easily performed by hand up to degree four and the
reader is strongly encouraged to do this.  The following degree five
result was obtained with a little bit of help from {\tt maple\/}:

\begin{align*}Z\left(T\left(m,p\right)\right)&=
         \exp_{\mbox{.}}\!\biggl(\!(m^2-1)(p^2-1)\biggl(\!-\tfrac{1 }{48}\vpic{w2A}\!
        +\tfrac{mp}{288}\vpic{tw2A}\!
 -\tfrac{9m^2p^2-m^2-p^2-1}{11520}\vpic{t2w2A}\\
&\phantom{{}={}}\phantom{\exp_\cdot\biggl(\ }    
+\tfrac{(m^2+1)(p^2+1)}{5760}\vpic{w4A}
 +\tfrac{mp(71m^2p^2-19m^2-19p^2-9)}{345600}\vpic{t3w2A}\\
&\phantom{{}={}}\phantom{\exp_\cdot\biggl(\ }    
+\tfrac{mp(m^2+1)(p^2+1)}{345600}\vpic{x3w2A}
- \tfrac{mp(2m^2p^2+m^2+p^2)}{17280}\vpic{tw4A}+\dots
\biggr)\biggr)
\end{align*}

Note that the coefficients of the round wheels,
$\vcenter{\epsffile{smallw2A.eps}}$  
and $\vcenter{\epsffile{smallw4A.eps}}$, are respectively
$-(m^2-1)(p^2-1)/48$ and $(m^4-1)(p^4-1)/5760$; this will be
generalized in the next section.
First however, the
methodology of this section leads to the next theorem which I had not
been able to prove using more na\"\i ve methods --- for example the methods
discussed in \cite{Willerton:VIsasPolys}.
\begin{theorem}
If $v$ is a rational, type $n$, Vassiliev knot invariant then, as a
function of $m$ and $p$, 
$v(T(m,p))$ is a polynomial, symmetric in $m$ and $p$ and of degree at
most $n$ in each of $m$ and $p$.
\end{theorem}
\begin{proof}
The symmetry in $m$ and $p$ follows from the fact
that the $(m,p)$-torus
knot and the $(p,m)$-torus knots are of the same knot-type (see
\cite{BurdeZieschang}).  Because the Kontsevich integral is universal
for Vassiliev invariants it suffices to prove that in the Kontsevich
integral the
coefficient of a basis symmetrized diagram has the properties described
in the theorem.  Fix a basis of symmetrized diagrams of degree $n$.
Pick any element of this basis.  By the linearity of change of basis
from round diagrams to symmetrized diagrams, the coefficient
of this element in 
  $\Omega.
  \exp_{\mbox{.}}\!\left(\tfrac{p}{m}\tfrac{1}{2}\thetadiag\right)$
is a
polynomial in $\left( \tfrac{p}{m}\right)$ of degree at most $n$.  
The action of $\psi^m$
will just be to multiply by some power of $m$.  By the symmetry in $m$
and $p$ the result must be a symmetric polynomial of
degree at most $n$  in $m$ and $p$.
\end{proof}
\section{The Alexander polynomial, wheel coefficients and torus knots.}
From work on the Melvin-Morton Theorem it is known that the Alexander
polynomial of a knot corresponds to the coefficients of the wheels in
its Kontsevich integral.  This is used to identify the
coefficients of the wheels for the torus knots.
\subsection{Identifying the Alexander polynomial in the Kontsevich
integral.} 
For a knot $K$, let $A_K(t)$ be the Alexander polynomial normalized so
that $A_K(t)=A_K(-t)$ and $A_K(1)=1$.  Define the weight system $W_{\text{AC}}\colon \A \to
\Q[[h]]$ by first defining it on connected symmetrized diagrams as
$$W_{\text{AC}}(D) = 
 \left\{ \begin{array}{ll} 
           - 2 h^{2n} & \text{if\ } D= w_{2n},\\
           0 & \text{if $D$ has less than deg$D$ legs,}
 \end{array} \right.$$ 
then extend it linearly to the space of primitive elements and
multiplicatively, with 
respect to the connect sum product, to the
whole of $\A$ (or equivalently \cite{Kricker:AlexConway}, 
extend $W_{\text{AC}}$ multiplicatively
with respect to the 
$\amalg$ product).

The AC in $W_{\text{AC}}$ stands for Alexander-Conway, the reason
being given by the next theorem which follows from work on the
Melvin-Morton Theorem.
\begin{theorem}[\cite{BarGar:MelvinMorton,KrickerSpenceAitchison:Cabling}]
If $K$ is a knot then, as formal power series in $h$,
   $$A_K(e^h)= W_{\text{AC}}\circ Z(K).$$
\end{theorem}
\noindent This can be used to identify the coefficients of the wheels in
the Kontsevich integral in the following manner.  Take logarithms of
both sides, as $W_{\text{AC}}$ is multiplicative with respect to the
disjoint union product, one finds
$$\ln A_K(e^h)= W_{\text{AC}}\circ\ln_\amalg\circ Z(K).$$
From the definition of $W_{\text{AC}}$ and the fact that
$\ln_\amalg\circ Z(K)$ can be written as the sum of connected
symmetrized diagrams, one deduces:
\begin{theorem}
With respect to a polynomial basis of symmetrized diagrams, the
coefficient of $w_{2n}$ in $Z(K)$ is equal to $-\tfrac{1}{2}$ times
the coefficient of $h^{2n}$ in $\ln A_K(e^h)$.
\end{theorem}
Note that if a basis of connected round diagrams is taken which in
degree $2n$ consists of the round wheel together with diagrams with
fewer univalent vertices, then the coefficient of the round wheel with $n$
legs is the same as the coefficient of $w_{2n}$ in the above
theorem. This can now be applied to the case of torus knots.

\subsection{Wheels coefficients for torus knots.}

For the torus knot $T(m,p)$ it is known that the Alexander
polynomial is given by 
  $$A_{T(m,p)}(t)=\frac{(t^{mp/2}-t^{-mp/2})(t^{1/2}-t^{-1/2})}
                       {(t^{m/2}-t^{-m/2})(t^{p/2}-t^{-p/2})}.$$
Substituting $t=e^h$, rewriting the bracketted terms as $\sinh$s,
dividing the numerator and the denominator of the right hand side by $h^2mp$,
taking logarithms of both sides and recalling the definition of the
modified Bernoulli numbers $\{b_{2n}\}$ from Section~5, one obtains
$$\ln \left(A_{T(m,p)}(e^h)\right)
   =\textstyle\sum b_{2n}(m^{2n}-1)(p^{2n}-1)2h^{2n}.$$
It then follows from the theorem above
that the coefficient of $w_{2n}$, the wheel with $2n$
legs, in $Z(T(m,p))$ is precisely
$-(m^{2n}-1)(p^{2n}-1)b_{2n}$. 
This can be seen to be in accord with the calculation of Section~8.
\section*{Acknowledgements.}
Much of this work was done whilst I was supported by a Royal Society
Research Fellowship at the
University of Melbourne.

I would like to thank Sergei Chmutov and Andrew Kricker for useful
conversations, Jacob Mostovoy and Jens Lieberum for comments on
earlier versions of this paper, and Thang Le for sharing his
unpublished work with me.

\bibliographystyle{amsplain}
\bibliography{Vas}

\appendix
\section{Change of basis data.}
Here is the change of basis data, up to degree five, between round
diagrams and symmetrized diagrams, used in the calculation of
Section~8.  It was worked out by hand by expanding each symmetrized
diagram and then using the STU relations to rewrite in terms of the
round diagram basis.  The $\Lambda$-module structure certainly
simplified the calculations.  The calculations took more than a week
and the degree six result seems to be beyond the limit of my stamina.
I don't know how much could be gained by automating the procedure, but
I think there is still much structure to be exploited, for instance
there is the following.

Note how sparse these matrices are.  The presence of identity matrices
along the diagonal is explained by the fact mentioned above that a
round diagram can be written as the symmetrized diagram obtained by
removing the Wilson loop plus diagrams with fewer legs.  Le
\cite{Le:Unpublished} conjectured that the connect sum of $l$ round
diagrams can be written as a linear combination of symmetrized
diagrams with at least $l$ legs:  this would explain some of the
zeroes in the bottom left corners.  The bottom row could be calculated if
it was known, as was also suggested by Le, that $\exp_{\mbox{.}}\!(\thetadiag
/2)=\Omega\amalg\exp_\amalg(\,|\,/2)$.   Also some of the zeroes seem to
be reflecting 
$\Lambda$-module structure, but I can not yet make that more precise. 
\begin{xxalignat}{1}\intertext{Degree zero:}
\vpic{1A}&=1. \\
\intertext{Degree one:}
\vpic{IA}&=\vpic{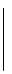}.\\
\intertext{Degree two:}
\begin{pmatrix}\vpic{w2A}\\ \vpic{I2A}\end{pmatrix}&=
\begin{pmatrix} \str 1&0\\ \str \tfrac{1}{6} &1 \end{pmatrix}
\begin{pmatrix}\vpic{w2B}\\ \vpic{IB}\vpic{IB}\end{pmatrix}.\\
\intertext{Degree three:}
 \begin{pmatrix}\vpic{tw2A}\\ \vpic{Iw2A}\\\vpic{I3A}\end{pmatrix}
&=  \begin{pmatrix} \str 1&0&0\\\str \tfrac{1}{3} &1&0\\ 
    \str 0&\tfrac{1}{2}&1 \end{pmatrix}
\begin{pmatrix}\vpic{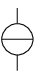}\\  \vpic{IB}\vpic{w2B}\\
       \vpic{IB}\vpic{IB} \vpic{IB}\end{pmatrix}.
\end{xxalignat}
Degree four:
\begin{xxalignat}{1}
\begin{pmatrix}\vpic{t2w2A}\\ \vpic{w4A}\\ \vpic{w2w2A}\\
       \vpic{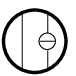}\\\vpic{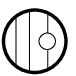}\\ \vpic{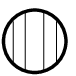}\end{pmatrix}
&= \begin{pmatrix} \str 1&0&0&0&0&0\\\str \tfrac{1}{3} &1&0&0&0&0\\ 
  \str \tfrac{2}{3}&0&1&0&0&0\\ \str \tfrac{1}{3}&0&0&1&0&0\\
  \str  0&0&\tfrac{1}{6}&\tfrac{2}{3}&1&0\\
  \str 0&-\tfrac{1}{15}&\tfrac{1}{12}&0&1&1
   \end{pmatrix}
\begin{pmatrix}\vpic{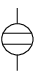}\\ \vpic{w4B}\\ \vpic{w2B}\vpic{w2B}\\
       \vpic{IB}\vpic{tw2B}\\\vpic{IB}\vpic{IB}\vpic{w2B}\\
       \vpic{IB}\vpic{IB} \vpic{IB} \vpic{IB}\end{pmatrix}.
\end{xxalignat}
Degree five:
\begin{xxalignat}{1}
\begin{pmatrix}\vpic{t3w2A}\\\vpic{x3w2A}\\ \vpic{tw4A}\\ \vpic{tw2w2A}\\
       \vpic{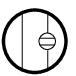}\\ \vpic{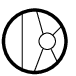}\\\vpic{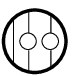}\\
        \vpic{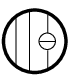}\\ \vpic{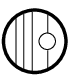}\\ \vpic{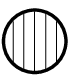}\end{pmatrix}
&= \begin{pmatrix} \str 1 & 0 & 0 & 0 & 0 & 0 & 0 & 0 & 0 & 0 \\
\str 0 & 1 & 0 & 0 & 0 & 0 & 0 & 0 & 0 & 0 \\
{\str \tfrac {1}{3}} & 0 & 1 & 0 & 0 & 0 & 0 & 0 & 0 & 0
 \\ \str
{\tfrac {2}{3}} & 0 & 0 & 1 & 0 & 0 & 0 & 0 & 0 & 0
 \\ \str
{\tfrac {1}{3}} & 0 & 0 & 0 & 1 & 0 & 0 & 0 & 0 & 0
 \\ \str
{\tfrac {2}{15}} & -{\tfrac {1}{30}} & 
{\tfrac {2}{3}} & 0 & {\tfrac {1}{3}}
 & 1 & 0 & 0 & 0 & 0 \\ \str
0 & 0 & 0 & {\tfrac {2}{3}} & {\tfrac {
2}{3}} & 0 & 1 & 0 & 0 & 0 \\\str 
0 & 0 & 0 & {\tfrac {1}{6}} & {\tfrac {
2}{3}} & 0 & 0 & 1 & 0 & 0 \\ \str
0 & 0 & -{\tfrac {2}{15}} & {\tfrac {1
}{6}} & 0 & 0 & {\tfrac {1}{2}} & 1 & 1 & 0 \\\str 
0 & 0 & 0 & 0 & 0 & -{\tfrac {1}{3}} & 
{\tfrac {5}{12}} & 0 & { \tfrac {5}{3}}
 & 1
  \end{pmatrix}
\begin{pmatrix}\vpic{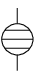}\\\vpic{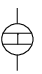}\\  \vpic{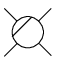}\\
\vpic{tw2B}\vpic{w2B}\\ 
       \vpic{IB}\vpic{t2w2B}\\ \vpic{IB}\vpic{w4B}\\
       \vpic{IB}\vpic{w2B}\vpic{w2B}\\
       \vpic{IB}\vpic{IB}\vpic{tw2B}\\\vpic{IB}\vpic{IB}\vpic{IB}\vpic{w2B}\\
       \vpic{IB}\vpic{IB} \vpic{IB} \vpic{IB}\vpic{IB}\end{pmatrix}.
\end{xxalignat}

\section{Comparison with Alvarez and Labastida's results}
In \cite{AlvarezLabastida:Torus}, Alvarez and Labastida calculate
Vassiliev invariants up to degree six for the torus knots by using the
values on the torus knots of the HOMFLY, Kauffman and Akutsu-Wadati
polynomials, together with a Chern-Simons approach to the universal
Vassiliev invariant.  Their method will not give all Vassiliev
invariants of torus knots because the invariants coming from Lie
algebras do not span the space of Vassiliev invariants
\cite{Vogel:AlgebraicStructures} (although they do so up to at
least degree twelve \cite{Kneissler:VasDegreeTwelve}).  However, up to
degree five, their results suggest that $\ln_{\mbox{.}}\! Z(T(m,p))$ might be
the following:
\begin{align*}
&(m^2-1)(p^2-1)\biggl(\tfrac{1}{6} \vpic{tree3A} +
         \tfrac{mp}{18} \vpic{tree4A}+
         \tfrac{9m^2p^2-m^2-p^2-1}{360} \vpic{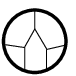}\\
    &\phantom{(m^2-1)(p^2-1)(}   +\tfrac{(m^2+1)(p^2+1)}{360}\vpic{w4A} 
  +  \tfrac{mp(69m^2p^2-21m^2-21p^2-11)}{5400}\vpic{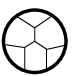}\\
  &\phantom{(m^2-1)(p^2-1)(}
 + \tfrac{mp(11m^2p^2+m^2+p^2-9)}{5400}\vpic{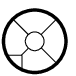}
   +\tfrac{mp(m^2+1)(p^2+1)}{900}\vpic{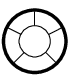}\biggr).
\end{align*}  
Most of the basis element for $\A$ which they use differ only by a
scalar from the basis elements used in the calculations above;
the only non-trivial change of basis relation is%
\footnote{One of the signs given in this expression in
\cite{ChmutovVarchenko} is incorrect.}
$$\vpic{w5A}=\tfrac{3}{2}\vpic{tw4A}-\tfrac{1}{6}\vpic{t3w2A}
-\tfrac{1}{12}\vpic{x3w2A}.$$ 
Rewriting, in the basis of Alvarez and Labastida, the expression
derived in Section~4, the following expression for $\ln_{\mbox{.}}\!
Z(T(m,p))$ up to degree five is obtained: 
\begin{align*}
&(m^2-1)(p^2-1)\biggl(
 -\tfrac{1}{24} \vpic{tree3A} +
         \tfrac{mp}{144} \vpic{tree4A}-
         \tfrac{9m^2p^2-m^2-p^2-1}{5760}
\vpic{tree5A}\\
 &\phantom{(m^2-1)(p^2-1)(} +\tfrac{(m^2+1)(p^2+1)}{5760}\vpic{w4A}
      +
   \tfrac{mp(69m^2p^2-21m^2-21p^2-11)}{172800}\vpic{tree6A}\\
   &\phantom{(m^2-1)(p^2-1)(}
 - \tfrac{mp(11m^2p^2+m^2+p^2-9)}{172800}\vpic{tailw4A}
  -\tfrac{mp(m^2+1)(p^2+1)}{28800}\vpic{w5A}\biggr).
\end{align*}
Note that this expression differs from that of Alvarez and Labastida
in the following two respects:  for a diagram of degree $n$, there is
a factor of 
$2^n$ difference --- this is probably explained taking a different norm
on the Lie algebras; and diagrams with an odd number of external
vertices differ in sign --- this is probably explained by an
orientation convention.

\end{document}